# The Generalized-Alpha-Beta-Skew-Normal Distribution: Properties and Applications

Sricharan Shah[1], Subrata Chakraborty[1], Partha Jyoti Hazarika[1] and M. MasoomAli[2]

[1]Department of Statistics, Dibrugarh University, Dibrugarh, Assam 786004 India

[2]Department of Mathematical Sciences, Ball State University, Muncie, IN 47306 USA

**Abstract**

In this paper we have introduced a generalized version of alpha beta skew normal distribution in the same line of Sharafi et al. (2017) and investigated some of its basic properties. The extensions of the proposed distribution have also been studied. The appropriateness of the proposed distribution has been tested by comparing the values of Akaike Information Criterion (AIC) and Bayesian Information Criterion (BIC) with the values of some other known related distributions for better model selection. Likelihood ratio testhas been used for discriminating between nested models.

**Math Subject Classification:** 60E05, 62H10, 62H12

1. **Introduction**

In many real life situations the data exhibit many modes as well as asymmetry, e.g., in the fields of demography, insurance, medical sciences, physics, etc. (for details see Hammel et al. (1983), Dimitrov et al. (1997), Sinha (2012), Chakraborty and Hazarika (2011), Chakraborty et al. (2015), among others). Azzalini (1985) first introduced the skew normal distribution and the density function of this distribution, denoted by $SN(\lambda)$, is given by

$$f_Z(z;\lambda) = 2\varphi(z)\Phi(\lambda z); \quad -\infty < z < \infty, -\infty < \lambda < \infty \quad (1)$$

where, $\varphi(.)$ is the density function of standard normal distribution, $\Phi(.)$ is the cumulative distribution function (cdf) of standard normal distribution and $\lambda$ is the asymmetry parameter.

The generalization of the skew normal distribution was proposed by Balakrishnan (2002) as a discussant in Arnold and Beaver (2002) and studied its properties and the density function is given by

$$f_Z(z;\lambda,n) = \varphi(z)[\Phi(\lambda z)]^n \big/ C_n(\lambda) \quad (2)$$

where, $n$, is a positive integer and $C_n(\lambda) = E[\Phi^n(\lambda U)]$, $U \sim N(0,1)$.

The general formula for the construction of skew-symmetric distributions was proposed by Huang and Chen (2007) starting from a symmetric (about 0) pdf $h(.)$ by introducing the concept of skew function $G(.)$, a Lebesgue measurable function such that, $0 \leq G(z) \leq 1$ and $G(z) + G(-z) = 1$, $z \in R$, almost everywhere. The density function of the same is given by

$$f(z) = 2h(z)G(z) \; ; z \in R. \tag{3}$$

Elal-Olivero (2010) introduced a new form of skew distribution which has both unimodal as well as bimodal behavior and is known as alpha skew normal distribution, denoted by $ASN(\alpha)$ and its density function is given by

$$f(z;\alpha) = \left(\frac{(1-\alpha z)^2 + 1}{2 + \alpha^2}\right)\varphi(z); z \in R. \tag{4}$$

Along the same line of the eqn. (4), Venegas et al. (2016) studied the logarithmic form of alpha-skew normal distribution and used for modeling chemical data. Sharafi et al. (2017) introduced a generalization of $ASN(\alpha)$ distribution. Shafiei et al. (2016) introduced a new family of skew distributions with more flexibility than the Azzalini (1985) and the Elal-Olivero (2010) distributions and is defined as follows.

A random variable $Z$ is said to be an alpha-beta skew normal distribution, denoted by $ABSN(\alpha, \beta)$ if its density function is given by

$$f(z;\alpha,\beta) = \left(\frac{(1-\alpha z - \beta z^3)^2 + 1}{2 + \alpha^2 + 15\beta^2 + 6\alpha\beta}\right)\varphi(z); z,\alpha,\beta \in R. \tag{5}$$

The reasons for considering this new distribution are as follows: First, this distribution includes four important classes of distribution which are normal, skew normal, alpha-beta skew normal, and the generalization of $ASN(\alpha)$ distribution. Second, this distribution is flexible enough to support both unimodality and bimodality (or multimodality) behaviors and has at most four modes. Third, for some real life data, this new proposed distribution provides better fitting as compared to the other distributions considered. The article is summarized as follows: In Section 2 we define the proposed distribution and some of its basic properties are investigated. In Section 3 we study some of its important distributional properties. Some extension of this distribution is discussed in Section 4. A stochastic representation and simulation procedure is given in Section 5. The parameter estimation and the real life applications of the proposed distribution are provided in Section 6. Lastly, conclusions are given in Section 7.

## 2. The Generalized Alpha-Beta Skew Normal Distribution

In this section we introduce a new generalized form of alpha-beta skew normal distribution and investigate some of its basic properties.

**Definition 1:** If a random variable $Z$ has a density function

$$f(z;\alpha,\beta,\lambda) = \frac{[(1-\alpha z - \beta z^3)^2 + 1]}{C(\alpha,\beta,\lambda)} \varphi(z)\Phi(\lambda z); \; z \in R \qquad (6)$$

where $C(\alpha,\beta,\lambda) = 1 + 3\alpha\beta - \alpha b\delta - \beta b\delta \frac{3+2\lambda^2}{1+\lambda^2} + \frac{\alpha^2}{2} + \frac{15\beta^2}{2}$, then it is said to be generalized alpha-beta skew normal distribution with skewness parameters $(\alpha,\beta,\lambda) \in R$. We denote it as $GABSN(\alpha,\beta,\lambda)$. The normalizing constant $C(\alpha,\beta,\lambda)$ is calculated as follows:

$$C(\alpha,\beta,\lambda) = \int_{-\infty}^{\infty} [(1-\alpha z - \beta z^3)^2 + 1]\varphi(z)\Phi(\lambda z)\,dz$$

$$= \int_{-\infty}^{\infty} (2 - 2\alpha z + \alpha^2 z^2 - 2\beta z^3 + 2\alpha\beta z^4 + \beta^2 z^6)\varphi(z)\Phi(\lambda z)\,dz$$

$$= \int_{-\infty}^{\infty} \varphi(z;\lambda)\,dz - \alpha E(Z_\lambda) + \frac{\alpha^2}{2}E(Z_\lambda^2) - \beta E(Z_\lambda^3) + \alpha\beta E(Z_\lambda^4) + \frac{\beta^2}{2}E(Z_\lambda^6)$$

$$= 1 + \frac{\alpha^2}{2} + 3\alpha\beta + \frac{15\beta^2}{2} - \alpha\sqrt{\frac{2}{\pi}}\frac{\lambda}{\sqrt{1+\lambda^2}} - \beta\sqrt{\frac{2}{\pi}}\frac{\lambda}{\sqrt{1+\lambda^2}}\frac{3+2\lambda^2}{1+\lambda^2}$$

$$= 1 + 3\alpha\beta - \alpha b\delta - \beta b\delta\frac{3+2\lambda^2}{1+\lambda^2} + \frac{\alpha^2}{2} + \frac{15\beta^2}{2},$$

where $b = \sqrt{\frac{2}{\pi}}, \delta = \frac{\lambda}{\sqrt{1+\lambda^2}}$ and $\varphi(z;\lambda)$ is the density function of $Z_\lambda \sim SN(\lambda)$.

### 2.1. Properties of $GABSN(\alpha,\beta,\lambda)$:

- If $\beta = 0$, then we get the generalized $ASN(\alpha)$ distribution of Sharafi et al. (2017) given by

$$f(z;\alpha,\lambda) = \frac{[(1-\alpha z)^2 + 1]}{1 - \alpha b\delta + \frac{\alpha^2}{2}}\varphi(z)\Phi(\lambda z)$$

- If $\alpha = 0$, then we get $f(z;\beta,\lambda) = \frac{[(1-\beta z^3)^2 + 1]}{1 - \beta b\delta\frac{3+2\lambda^2}{1+\lambda^2} + \frac{15\beta^2}{2}}\varphi(z)\Phi(\lambda z)$

  This is known as generalized beta skew normal $GBSN(\beta,\lambda)$ distribution.

- If $\lambda = 0$, then we get the $ABSN(\alpha,\beta)$ distribution of Shafiei et al. (2016) given by

$$f(z;\alpha,\beta) = \frac{[(1-\alpha z - \beta z^3)^2 + 1]}{2 + \alpha^2 + 6\alpha\beta + 15\beta^2}\varphi(z)$$

- If $\alpha = \beta = 0$, then we get the $SN(\lambda)$ distribution of Azzalini (1985) and is given by

$$f(z;\lambda) = 2\varphi(z)\Phi(\lambda z)$$

- If $\alpha = \beta = \lambda = 0$, then we get the standard normal distribution and is given by

$$f(z) = \varphi(z)$$

- If $\alpha \to \pm\infty$, then $f(z;\lambda) = 2z^2\varphi(z)\Phi(\lambda z)$, where, $2z^2\varphi(z)$ is the pdf of $BN(2)$ (see Hazarika et al. (2019)). Therefore, as $\alpha \to \pm\infty$, $GABSN(\alpha,\beta,\lambda) \to GBN(2)$

- If $\beta \to \pm\infty$, then $GABSN(\alpha,\beta,\lambda) \to GBN(6)$. the pdf of $GBN(6)$ is

$$f(z;\lambda) = \frac{2z^6}{15}\varphi(z)\Phi(\lambda z)$$

- For fixed $\alpha$ and $\beta$, if $\lambda \to +\infty$ then

$$f(z;\alpha,\beta) = \frac{[(1-\alpha z - \beta z^3)^2 + 1]}{1 - \alpha b - 2\beta b + 3\alpha\beta + \frac{\alpha^2}{2} + \frac{15\beta^2}{2}}\varphi(z)I(z>0)$$

and if $\lambda \to -\infty$ then

$$f(z;\alpha,\beta) = \frac{[(1-\alpha z - \beta z^3)^2 + 1]}{1 + \alpha b + 2\beta b + 3\alpha\beta + \frac{\alpha^2}{2} + \frac{15\beta^2}{2}}\varphi(z)I(z<0)$$

- If $Z \sim GABSN(\alpha,\beta,\lambda)$, then $-Z \sim -GABSN(\alpha,\beta,\lambda)$.

## 2.2. Plots of the density function

The density function of $GABSN(\alpha,\beta,\lambda)$ distribution for different choices of the parameters $\alpha, \beta$ and $\lambda$ are plotted in Figure 1.

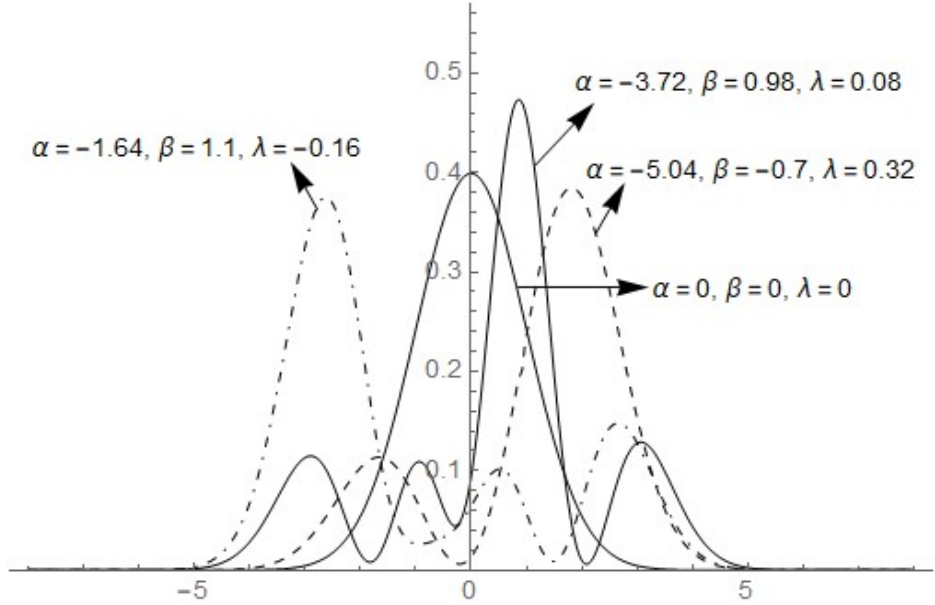

**Figure 1:** Plots of the density function of $GABSN(\alpha, \beta, \lambda)$

## 2.3. Mode of $GABSN(\alpha, \beta, \lambda)$:

**Theorem 1:** The $GABSN(\alpha, \beta, \lambda)$ distribution has at most four modes.

*Proof:* Let $Z \sim GABSN(\alpha, \beta, \lambda)$ distribution then by eqn. (5) and eqn. (6) we have,

$$f(z;\alpha,\beta,\lambda) = \frac{2+\alpha^2+15\beta^2+6\alpha\beta}{C(\alpha,\beta,\lambda)} f(z;\alpha,\beta)\Phi(\lambda z) . \quad (7)$$

Then by differentiating we get,

$$f'(z;\alpha,\beta,\lambda) = \frac{2+\alpha^2+15\beta^2+6\alpha\beta}{C(\alpha,\beta,\lambda)} \left[ f'(z;\alpha,\beta)\Phi(\lambda z) + \lambda f(z;\alpha,\beta)\varphi(\lambda z) \right] . \quad (8)$$

To prove that the eqn. (7) has at most four modes, we have to show that the eqn. (8) has one or seven roots. For this, we apply a graphical approach, and so we have

$$f'(z;\alpha,\beta,\lambda) = F_1(z) - F_2(z)$$

where

$$F_1(z) = \frac{2+\alpha^2+15\beta^2+6\alpha\beta}{C(\alpha,\beta,\lambda)} f'(z;\alpha,\beta)\Phi(\lambda z)$$

$$F_2(z) = -\frac{2+\alpha^2+15\beta^2+6\alpha\beta}{C(\alpha,\beta,\lambda)} \lambda f(z;\alpha,\beta)\varphi(\lambda z) .$$

Then, set $f'(z;\alpha,\beta,\lambda) = 0$, and hence $F_1(z) = F_2(z)$. (9)

Since the eqn. (6) vanishes for (-5, 5), we draw the curves of $C1: y = F_1(z)$ for $\alpha = -3.72, \beta = 1.08, \lambda = 0.3$ and $C2: y = F_2(z)$ for $\alpha = -1.53, \beta = 0.1, \lambda = 3$ in Figure 2.

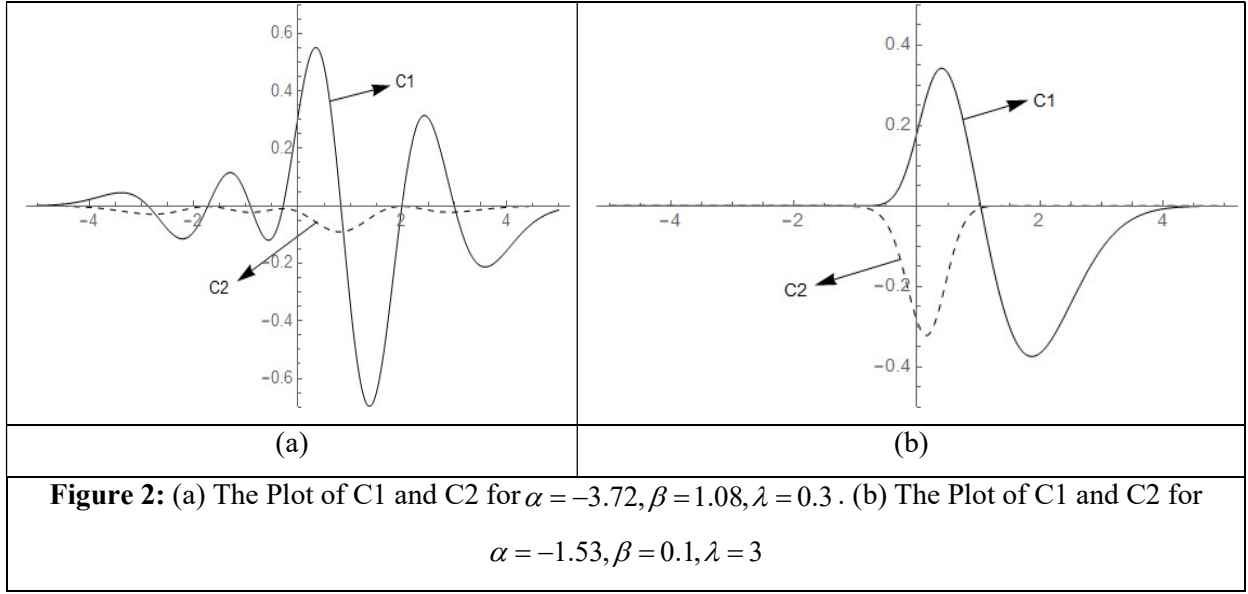

| (a) | (b) |

**Figure 2:** (a) The Plot of C1 and C2 for $\alpha = -3.72, \beta = 1.08, \lambda = 0.3$. (b) The Plot of C1 and C2 for $\alpha = -1.53, \beta = 0.1, \lambda = 3$

From the figure 2, it is shown that the two curves have atleast one and at most seven intersection points and the values of $z$ of these points are the roots of eqn. (9).

Since $\lim_{z \to \pm\infty} f(z;\alpha,\beta,\lambda) = 0$, then if eqn.(9) has one root it should be the mode of eqn. (7) and if eqn.(9) has seven roots, then eqn. (7) should have four modes and hence the $GABSN(\alpha,\beta,\lambda)$ distribution has at least one and at most four modes.

## 3. Distributional Properties

### 3.1. Cumulative Distribution Function

**Theorem 1:** The cdf of $GABSN(\alpha,\beta,\lambda)$ distribution is given by

$$F_Z(z) = \frac{1}{2(1+\lambda^2)^{5/2}}(b\delta(\alpha^2(1+\lambda^2)^2 - 2\alpha\beta(1+\lambda^2)(5+z^2+(3+z^2)\lambda^2) - \beta(-2z(1+\lambda^2)^2 + z^4\beta(1+\lambda^2)^2 +$$

$$z^2\beta(1+\lambda^2)(9+5\lambda^2) + \beta(33+40\lambda^2+15\lambda^4)))\varphi(z\sqrt{1+\lambda^2}) - \frac{1}{\sqrt{\pi}}(1+\lambda^2)^2(\sqrt{2}(\alpha+2\beta)\lambda\, Erf(z\sqrt{1+\lambda^2}/\sqrt{2})$$

$$+ \sqrt{\pi}(2(\alpha^2 z + 2\alpha(-1+z(3+\lambda^2)\beta) + \beta(-2(2+\lambda^2)+z(15+5z^2+z^4)\beta))\sqrt{1+\lambda^2}\varphi(z)\Phi(\lambda z) + 2b\beta\delta\,\Phi(z\sqrt{1+\lambda^2})$$

$$- 3(1+2\alpha\beta+5\beta^2)\sqrt{1+\lambda^2}\,\Phi(z;\lambda))))$$

(10)

*Proof:* $F_Z(z) = P(Z \le z) = \int_{-\infty}^{z} \frac{(1-\alpha t - \beta t^3)^2 + 1}{C(\alpha,\beta,\lambda)} \varphi(t)\Phi(\lambda t)dt$

$$= \frac{1}{C(\alpha,\beta,\lambda)}\int_{-\infty}^{z}(2 - 2\alpha t + \alpha^2 t^2 - 2\beta t^3 + 2\alpha\beta t^4 + \beta^2 t^6)\varphi(t)\Phi(\lambda t)dt$$

$$= \frac{1}{C(\alpha,\beta,\lambda)}\left[\int_{-\infty}^{z} 2\varphi(t)\Phi(\lambda t)dt - 2\alpha\int_{-\infty}^{z} t\varphi(t)\Phi(\lambda t)dt + \alpha^2\int_{-\infty}^{z} t^2\varphi(t)\Phi(\lambda t)dt - 2\beta\int_{-\infty}^{z} t^3\varphi(t)\Phi(\lambda t)dt + \right.$$

$$\left. 2\alpha\beta\int_{-\infty}^{z} t^4\varphi(t)\Phi(\lambda t)dt + \beta^2\int_{-\infty}^{z} t^6\varphi(t)\Phi(\lambda t)dt \right]$$

(11)

Since, we have

$$\int_{-\infty}^{z} 2\varphi(t)\Phi(\lambda t)dt = \Phi(z;\lambda)$$

$$\int_{-\infty}^{z} t\varphi(t)\Phi(\lambda t)dt = -\varphi(z)\Phi(\lambda z) + \frac{\lambda\, Erf\left((z\sqrt{1+\lambda^2})/\sqrt{2}\right)}{2\sqrt{2\pi}\sqrt{1+\lambda^2}}$$

$$\int_{-\infty}^{z} t^2\varphi(t)\Phi(\lambda t)dt = -z\varphi(z)\Phi(\lambda z) + \frac{b\delta\,\varphi\left(z\sqrt{1+\lambda^2}\right)}{2\sqrt{1+\lambda^2}}$$

$$\int_{-\infty}^{z} t^3\varphi(t)\Phi(\lambda t)dt = -(2+z^2)\varphi(z)\Phi(\lambda z) + \frac{\lambda\, Erf\left((z\sqrt{1+\lambda^2})/\sqrt{2}\right)}{\sqrt{2\pi}\sqrt{1+\lambda^2}} - z\frac{b\delta\,\varphi\left(z\sqrt{1+\lambda^2}\right)}{2\sqrt{1+\lambda^2}} + \frac{b\delta\,\Phi\left(z\sqrt{1+\lambda^2}\right)}{2\sqrt{1+\lambda^2}}$$

$$\int_{-\infty}^{z} t^4\varphi(t)\Phi(\lambda t)dt = -z(3+z^2)\varphi(z)\Phi(\lambda z) - \frac{b\delta\,\varphi\left(z\sqrt{1+\lambda^2}\right)(5+z^2+\lambda^2(3+z^2))}{2\left(\sqrt{1+\lambda^2}\right)^3} + \frac{3}{2}\Phi(z;\lambda)$$

$$\int_{-\infty}^{z} t^6\varphi(t)\Phi(\lambda t)dt = -z(15+5z^2+z^4)\varphi(z)\Phi(\lambda z) + \frac{b\delta\varphi\left(z\sqrt{1+\lambda^2}\right)\left(-33-40\lambda^2-15\lambda^4-z^4(1+\lambda^2)^2-z^2(9+14\lambda^2+5\lambda^4)\right)}{2\left(\sqrt{1+\lambda^2}\right)^5} + \frac{15}{2}\Phi(z;\lambda)$$

Putting the above results in eqn. (11), we get the desired result in eqn. (10).

## 3.2. Moments

**Theorem 3:** The $k^{th}$ order moment of $GABSN(\alpha,\beta,\lambda)$ distribution is given by

$$E(Z_\lambda^k) = \frac{1}{C(\alpha,\beta,\lambda)}\left[E(Z_\lambda^k) - \alpha\, E(Z_\lambda^{k+1}) + \frac{\alpha^2}{2}E(Z_\lambda^{k+2}) - \beta\, E(Z_\lambda^{k+3}) + \alpha\beta\, E(Z_\lambda^{k+4}) + \frac{\beta^2}{2}E(Z_\lambda^{k+6})\right]$$

where, $E(Z_\lambda^k)$ is the $k^{th}$ moment of $Z_\lambda \sim SN(\lambda)$. (12)

*Proof:* $E(Z^k) = \int_{-\infty}^{\infty} z^k \frac{(1-\alpha z - \beta z^3)^2 + 1}{C(\alpha,\beta,\lambda)} \varphi(z)\Phi(\lambda z)\, dz$

$$= \frac{1}{C(\alpha,\beta,\lambda)}\left[\int_{-\infty}^{\infty} 2z^k\varphi(z)\Phi(\lambda z)dz - 2\alpha\int_{-\infty}^{\infty} z^{k+1}\varphi(z)\Phi(\lambda z)dz + \alpha^2\int_{-\infty}^{\infty} z^{k+2}\varphi(z)\Phi(\lambda z)dz\right.$$

$$\left. - 2\beta\int_{-\infty}^{\infty} z^{k+3}\varphi(z)\Phi(\lambda z)dz + 2\alpha\beta\int_{-\infty}^{\infty} z^{k+4}\varphi(z)\Phi(\lambda z)dz + \beta^2\int_{-\infty}^{\infty} z^{k+6}\varphi(z)\Phi(\lambda z)dz\right]$$

$$= \frac{1}{C(\alpha,\beta,\lambda)}\left[E(Z_\lambda^k) - 2\alpha\, E(Z_\lambda^{k+1}) + \alpha^2\, E(Z_\lambda^{k+2}) - 2\beta\, E(Z_\lambda^{k+3}) + 2\alpha\beta\, E(Z_\lambda^{k+4}) + \beta^2\, E(Z_\lambda^{k+5})\right]$$

The moments of $Z \sim GABSN(\alpha,\beta,\lambda)$ can be obtained by applying the moments of $Z_\lambda \sim SN(\lambda)$ (Henze, 1986). Thus, we obtained the following results for $k = 1, 2, 3, 4$:

$$E(X) = \frac{1}{C(\alpha,\beta,\lambda)}\left[-\alpha - 3\beta + b\delta + \frac{b\delta\alpha^2(3+2\lambda^2)}{2(1+\lambda^2)} + \frac{b\delta\alpha\beta(c_1)}{(1+\lambda^2)^2} + \frac{3b\delta\beta^2(c_2)}{2(1+\lambda^2)^3}\right]$$

$$E(X^2) = \frac{1}{C(\alpha,\beta,\lambda)}\left[1 + 15\alpha\beta + \frac{3\alpha^2}{2} + \frac{105\beta^2}{2} - \frac{b\delta\alpha(3+2\lambda^2)}{(1+\lambda^2)} - \frac{b\delta\beta(c_1)}{(1+\lambda^2)^2}\right]$$

$$E(X^3) = \frac{1}{C(\alpha,\beta,\lambda)}\left[-3\alpha - 15\beta + \frac{b\delta(3+2\lambda^2)}{(1+\lambda^2)} + \frac{b\delta\alpha^2(c_1)}{2(1+\lambda^2)^2} + \frac{3b\delta\alpha\beta(c_2)}{(1+\lambda^2)^3} + \frac{3b\delta\beta^2(c_3)}{2(1+\lambda^2)^4}\right]$$

$$E(X^4) = \frac{1}{C(\alpha,\beta,\lambda)}\left[3 + 105\alpha\beta + \frac{15\alpha^2}{2} + \frac{945\beta^2}{2} - \frac{b\delta\alpha(c_1)}{(1+\lambda^2)^2} - \frac{3b\delta\beta(c_2)}{(1+\lambda^2)^3}\right]$$

$$Var(X) = \frac{1}{[C(\alpha,\beta,\lambda)]^2}\left[C(\alpha,\beta,\lambda)\left(1 + 15\alpha\beta + \frac{3\alpha^2}{2} + \frac{105\beta^2}{2} - \frac{b\delta\alpha(3+2\lambda^2)}{(1+\lambda^2)} - \frac{b\delta\beta(c_1)}{(1+\lambda^2)^2}\right) - \right.$$
$$\left.\left(-\alpha - 3\beta + b\delta + \frac{b\delta\alpha^2(3+2\lambda^2)}{2(1+\lambda^2)} + \frac{b\delta\alpha\beta(c_1)}{(1+\lambda^2)^2} + \frac{3b\delta\beta^2(c_2)}{2(1+\lambda^2)^3}\right)^2\right]$$

where, $c_1 = 15 + 20\lambda^2 + 8\lambda^4$; $c_2 = 35 + 70\lambda^2 + 56\lambda^4 + 16\lambda^6$; $c_3 = 315 + 8\lambda^2(105 + 126\lambda^2 + 72\lambda^4 + 16\lambda^6)$.

**Remark 1:** Using numerical optimization of $E(Z)$ and $Var(Z)$ with respect to $\alpha, \beta$ and $\lambda$, the following bounds for mean and variance can be obtained as $-2.75863 \leq E(Z) \leq 2.75863$ and $0.31173 \leq Var(Z) \leq 8.16228$

**Remark 2:** By taking the limit $\alpha \to \pm\infty$ in the moments of $GABSN(\alpha,\beta,\lambda)$ distribution, the moments of $GBN(2)$ distribution can be obtained as

$$E(Z) \to \frac{b\delta(3+2\lambda^2)}{(1+\lambda^2)}; \quad Var(Z) \to \frac{3\pi(1+\lambda^2)^3 - 2\lambda^2(3+2\lambda^2)^2}{\pi(1+\lambda^2)^3}.$$

**Remark 3:** By taking the limit $\beta \to \pm\infty$ in the moments of $GABSN(\alpha,\beta,\lambda)$ distribution, the moments of $GBN(6)$ distribution can be obtained easily as

$$E(Z) \to \frac{b\delta(35 + 70\lambda^2 + 56\lambda^4 + 16\lambda^6)}{5(1+\lambda^2)^3}; \quad Var(Z) = \frac{175\pi(1+\lambda^2)^7 - 2\lambda^2(35 + 70\lambda^2 + 56\lambda^4 + 16\lambda^6)^2}{25\pi(1+\lambda^2)^7}.$$

**Remark 4:** By taking the limit $\lambda \to +\infty$ or $-\infty$ in the moments of $GABSN(\alpha,\beta,\lambda)$ distribution, we get when, $\lambda \to +\infty$

$$E(Z) = \frac{b(-2 - 2\alpha^2 + \alpha(\sqrt{2\pi} - 16\beta) + 3\sqrt{2\pi}\beta - 48\beta^2)}{-2 + 2b\alpha - \alpha^2 + 4b\beta - 6\alpha\beta - 15\beta^2};$$

$$Var(Z) = -\frac{2(2 - \sqrt{2\pi}\alpha + 2\alpha^2 - 3\sqrt{2\pi}\beta + 16\alpha\beta + 48\beta^2)^2}{\pi(2 - 2b\alpha + \alpha^2 - 4b\beta + 6\alpha\beta + 15\beta^2)^2} + \frac{2 - 4b\alpha + 3\alpha^2 - 16b\beta + 30\alpha\beta + 105\beta^2}{2 - 2b\alpha + \alpha^2 - 4b\beta + 6\alpha\beta + 15\beta^2}.$$

And when $\lambda \to -\infty$ then

$$E(Z) = \frac{b(2 + 2\alpha^2 + \alpha(\sqrt{2\pi} + 16\beta) + 3\sqrt{2\pi}\beta + 48\beta^2)}{2 + 2b\alpha + \alpha^2 + 4b\beta + 6\alpha\beta + 15\beta^2},$$

$$Var(Z) = -\frac{2(2 + \sqrt{2\pi}\alpha + 2\alpha^2 + 3\sqrt{2\pi}\beta + 16\alpha\beta + 48\beta^2)^2}{\pi(2 + 2b\alpha + \alpha^2 + 4b\beta + 6\alpha\beta + 15\beta^2)^2} + \frac{2 + 4b\alpha + 3\alpha^2 + 16b\beta + 30\alpha\beta + 105\beta^2}{2 + 2b\alpha + \alpha^2 + 4b\beta + 6\alpha\beta + 15\beta^2}.$$

### 3.3. Skewness and Kurtosis

**Remark 5:** The skewness and kurtosis of $GABSN(\alpha, \beta, \lambda)$ distribution is obtained respectively, by using the formulae

$$\beta_1 = \frac{(E(Z^3) - 3E(Z^2)E(Z) + 2[E(Z)]^3)^2}{(E(Z^2) - [E(Z)]^2)^3} \text{ and } \beta_2 = \frac{E(Z^4) - 4E(Z^3)E(Z) + 6E(Z^2)[E(Z)]^2 - 3[E(Z)]^4}{(E(Z^2) - [E(Z)]^2)^2}$$

which cannot be expressed conveniently as the expression is very vast. The following bounds for skewness and kurtosis can be calculated by numerically optimizing $\beta_1$ and $\beta_2$ with respect to $\alpha, \beta$ and $\lambda$ as $0 \leq \beta_1 \leq 6.70451$ and $1.22732 \leq \beta_2 \leq 14.1965$.

**Remark 6:** The skewness and kurtosis of $GBN(2)$ distribution can be obtained easily by taking limit $\alpha \to \pm\infty$ in the results of $GABSN(\alpha, \beta, \lambda)$ distribution as

$$\beta_1 = \frac{2(-4\lambda^3(3 + 2\lambda^2)^3 + \pi\lambda(1 + \lambda^2)^2(12 + 25\lambda^2 + 10\lambda^4))^2}{(3\pi(1 + \lambda^2)^3 - 2\lambda^2(3 + 2\lambda^2)^2)^3}$$

$$\beta_2 = \frac{15\pi^2(1 + \lambda^2)^6 - 12\lambda^4(3 + 2\lambda^2)^4 + 4\pi(\lambda + \lambda^3)^2(-9 + 9\lambda^2 + 16\lambda^4 + 4\lambda^6)}{(3\pi(1 + \lambda^2)^3 - 2\lambda^2(3 + 2\lambda^2)^2)^2}.$$

**Remark 7:** The skewness and kurtosis of $GBN(6)$ distribution can be obtained easily by taking limit $\beta \to \pm\infty$ in the results of $GABSN(\alpha, \beta, \lambda)$ distribution as

$$\beta_1 = \frac{2\lambda^2(4\lambda^2(c_2)^3 - 25\pi(1 + \lambda^2)^6(c_4))^2}{15625\pi^3(1 + \lambda^2)^{21}(c_6)^3}, \quad \beta_2 = \frac{39375\pi^2(1 + \lambda^2)^{14} - 12\lambda^4(c_2)^4 + 500\pi\lambda^2(1 + \lambda^2)^6 c_2 c_5}{625\pi^2(1 + \lambda^2)^{14}(c_6)^2},$$

where, $c_4 = 420 + 1365\lambda^2 + 1638\lambda^4 + 936\lambda^6 + 208\lambda^8$, $c_5 = 21 + 105\lambda^2 + 126\lambda^4 + 72\lambda^6 + 16\lambda^8$, and

$c_6 = 7 - \frac{2\lambda^2(35 + 70\lambda^2 + 56\lambda^4 + 16\lambda^6)^2}{25\pi(1 + \lambda^2)^7}.$

**Note:** By taking the limit $\lambda \to +\infty \text{ or } -\infty$ in the results of $GABSN(\alpha, \beta, \lambda)$ distribution, we obtained the results but we are unable to express the expression conveniently as the expression is very vast.

### 4. The Stochastic Representation for $GABSN(\alpha, \beta, \lambda)$ distribution

**Theorem 3:** The conditional distribution of $W | \{\lambda W > X\}$ follows $GABSN(\alpha, \beta, \lambda)$ distribution, if $W \sim ABSN(\alpha, \beta)$ and $X \sim N(0,1)$, and are independent.

*Proof*: Assume $Z = W \mid \{\lambda W > X\}$. Then, we can have

$$P(Z \leq z) = P(W \leq z \mid \lambda W > X) = \frac{P(W \leq z, \lambda W > X)}{P(\lambda W > X)}$$

Using eqn. (5), we get, $P(W \leq z, \lambda W > X) = \int_{-\infty}^{z} \frac{(1-\alpha u - \beta u^3)^2 + 1}{2 + \alpha^2 + 15\beta^2 + 6\alpha\beta} \varphi(u) \Phi(\lambda u) du$

and $P(\lambda W > X) = \int_{-\infty}^{\infty} \frac{(1-\alpha u - \beta u^3)^2 + 1}{2 + \alpha^2 + 15\beta^2 + 6\alpha\beta} \varphi(u) \Phi(\lambda u) du = \frac{1}{2 + \alpha^2 + 15\beta^2 + 6\alpha\beta} C(\alpha, \beta, \lambda)$

Therefore, $P(Z \leq z) = \dfrac{\int_{-\infty}^{z} \frac{(1-\alpha u - \beta u^3)^2 + 1}{2 + \alpha^2 + 15\beta^2 + 6\alpha\beta} \varphi(u) \Phi(\lambda u) du}{\frac{1}{2 + \alpha^2 + 15\beta^2 + 6\alpha\beta} C(\alpha, \beta, \lambda)} = \int_{-\infty}^{z} \frac{(1-\alpha u - \beta u^3)^2 + 1}{C(\alpha, \beta, \lambda)} \varphi(u) \Phi(\lambda u) du$

and the density function of the conditional distribution of $Z = W \mid \{\lambda W > X\}$ is

$$f_Z(z) = \frac{(1-\alpha z - \beta z^3)^2 + 1}{C(\alpha, \beta, \lambda)} \varphi(z) \Phi(\lambda z) \sim GAbSN(\alpha, \beta, \lambda)$$

Or, we write $W \mid \{\lambda W > X\} \sim GAbSN(\alpha, \beta, \lambda)$.

**Algorithm**

We can simulate the data for the proposed distribution with the help of the above stochastic representation by applying the acceptance-rejection technique using the following steps:

a) Generate $X \sim N(0,1)$ and $W \sim ABSN(\alpha, \beta)$ (Sharafi et al., 2017) independently.

b) If $Z = W \mid \{\lambda W > X\}$ put $Z = W$. Otherwise, go to step a).

## 5. The Extension of Generalized Alpha Beta Skew Normal Distribution

### 5.1. Log-Generalized Alpha Beta Skew Normal Distribution

In this section, using the idea of (Venegas et al., 2016), we present the definition and some simple properties of log-generalized alpha beta skew normal distribution. Let $Z = e^Y$, then $Y = Log(Z)$, therefore, the density function of $Z$ is defined as follows:

**Definition 3:** If the random variable $Z$ has the density function given by

$$f(z; \alpha, \beta, \lambda) = \frac{[(1-\alpha y - \beta y^3)^2 + 1]}{C(\alpha, \beta, \lambda) z} \varphi(y) \Phi(\lambda y); \quad z > 0, \tag{13}$$

then we say that $Z$ is distributed according to the log-generalized alpha beta skew normal distribution with parameter $\alpha, \beta, \lambda \in R$, where, $y = Log(z)$ and $\varphi(y)$ is the density function of the standard log-normal distribution. We denote it by $LGABSN(\alpha, \beta, \lambda)$.

**Properties of $LGABSN(\alpha, \beta, \lambda)$:**

- If $\beta = 0$, then we get $f(z;\alpha,\lambda) = \dfrac{[(1-\alpha y)^2 + 1]}{\left(1 - \alpha b \delta + \dfrac{\alpha^2}{2}\right)z}\varphi(y)\Phi(\lambda y)$.

  This is known as log-generalized $ASN(\alpha)$ distribution.

- If $\alpha = 0$, then we get $f(z;\beta,\lambda) = \dfrac{[(1-\beta y^3)^2 + 1]}{\left(1 - \beta b \delta \dfrac{3+2\lambda^2}{1+\lambda^2} + \dfrac{15\beta^2}{2}\right)z}\varphi(y)\Phi(\lambda y)$.

  This is known as log-generalized beta skew normal $LGBSN(\beta,\lambda)$ distribution.

- If $\lambda = 0$, then we get $f(z;\alpha,\beta) = \dfrac{[(1-\alpha y - \beta y^3)^2 + 1]}{(2 + \alpha^2 + 6\alpha\beta + 15\beta^2)z}\varphi(y)$.

  This is known as log-generalized $ABSN(\alpha,\beta)$ distribution.

- If $\alpha = \beta = 0$, then we get $f(z;\lambda) = \dfrac{2\varphi(y)\Phi(\lambda y)}{z}$.

  This is known as log $SN(\lambda)$ distribution.

- If $\alpha = \beta = \lambda = 0$, then we get the standard log-normal distribution and is given by

$$f(z) = \dfrac{\varphi(y)}{z}$$

- If $\alpha \to \pm\infty$, then we get the log-generalized bimodal normal ($LGBN(2)$) distribution given by

$$f(z;\lambda) = \dfrac{y^2}{z}\varphi(y)\Phi(\lambda y)$$

- If $\beta \to \pm\infty$, then we get the log-generalized bimodal normal ($LGBN(6)$) distribution given by

$$f(z;\lambda) = \dfrac{y^6}{15z}\varphi(y)\Phi(\lambda y)$$

- If $Z \sim LGABSN(\alpha,\beta,\lambda)$, then $-Z \sim -LGABSN(\alpha,\beta,\lambda)$.

### 5.2. Location Scale Extension

The location and scale extension of $GABSN(\alpha,\beta,\lambda)$ distribution is as follows. If $Z \sim GABSN(\alpha,\beta,\lambda)$ then $Y = \mu + \sigma Z$ is said to be the location ($\mu$) and scale ($\sigma$) extension of $Z$ and has the density function is given by

$$f_Z(z;\alpha,\beta,\lambda,\mu,\sigma) = \dfrac{\left(\left[1 - \alpha\left(\dfrac{z-\mu}{\sigma}\right) - \beta\left(\dfrac{z-\mu}{\sigma}\right)^3\right]^2 + 1\right)}{C(\alpha,\beta,\lambda)}\varphi\left(\dfrac{z-\mu}{\sigma}\right)\Phi\left(\lambda\left(\dfrac{z-\mu}{\sigma}\right)\right), \qquad (14)$$

where, $(z, \alpha, \beta, \lambda, \mu) \in R$, and $\sigma > 0$. We denote it by $Y \sim GABSN(\alpha, \beta, \lambda, \mu, \sigma)$.

## 6. Parameter Estimation and Applications

### 6.1. Maximum Likelihood Estimation

Let $y_1, y_2, \ldots, y_n$ be a random sample from the distribution of the random variable $Y \sim GABSN(\alpha, \beta, \lambda, \mu, \sigma)$ so that the log-likelihood function for $\theta = (\alpha, \beta, \lambda, \mu, \sigma)$ is given by

$$l(\theta) = \sum_{i=1}^{n} \log\left[\left\{1 - \alpha\left(\frac{y_i - \mu}{\sigma}\right) - \beta\left(\frac{y_i - \mu}{\sigma}\right)^3\right\}^2 + 1\right] - n\log C(\alpha, \beta, \lambda) - n\log(\sigma) - \frac{n}{2}\log(2\pi) - \sum_{i=1}^{n} \frac{(y_i - \mu)}{\sigma} - 2\sum_{i=1}^{n} \log\left[\Phi\left(\frac{\lambda(y_i - \mu)}{\sigma}\right)\right]$$

(15)

Taking Partial derivatives of eqn. (15) w.r.t. the parameters, the following normal equations are obtained:

$$\frac{\partial l(\theta)}{\partial \mu} = \frac{n}{\sigma} + \sum_{i=1}^{n} \frac{2D\left(\frac{\alpha}{\sigma} + \frac{3\beta(y_i - \mu)^2}{\sigma^3}\right)}{1 + D^2} + \frac{2\lambda}{\sigma}\sum_{i=1}^{n} W\left(\frac{\lambda(y_i - \mu)}{\sigma}\right)$$

$$\frac{\partial l(\theta)}{\partial \sigma} = -\frac{n}{\sigma} + \sum_{i=1}^{n}\left(\frac{(y_i - \mu)}{\sigma^2}\right) + \sum_{i=1}^{n} \frac{2D\left(\frac{\alpha(y_i - \mu)}{\sigma^2} + \frac{3\beta(y_i - \mu)^3}{\sigma^4}\right)}{1 + D^2} + \frac{2\lambda}{\sigma^2}\sum_{i=1}^{n}(y_i - \mu)W\left(\frac{\lambda(y_i - \mu)}{\sigma}\right)$$

$$\frac{\partial l(\theta)}{\partial \alpha} = -\frac{n(\alpha + 3\beta - b\delta)}{C(\alpha, \beta, \lambda)} - \sum_{i=1}^{n} \frac{2D(y_i - \mu)}{\sigma(1 + D^2)}$$

$$\frac{\partial l(\theta)}{\partial \beta} = -\frac{n\left(3\alpha + 15\beta - \frac{b\delta(3 + 2\lambda^2)}{(1 + \lambda^2)}\right)}{C(\alpha, \beta, \lambda)} - \sum_{i=1}^{n} \frac{2D(y_i - \mu)^3}{\sigma^3(1 + D^2)}$$

$$\frac{\partial l(\theta)}{\partial \lambda} = \frac{n\left(\frac{b(\alpha + 3\beta + \alpha\lambda^2)}{(1 + \lambda^2)^{5/2}}\right)}{C(\alpha, \beta, \lambda)} - \frac{2}{\sigma}\sum_{i=1}^{n}(y_i - \mu)W\left(\frac{\lambda(y_i - \mu)}{\sigma}\right)$$

where, $D = 1 - \frac{\alpha(y_i - \mu)}{\sigma} - \frac{\beta(y_i - \mu)^3}{\sigma^3}$ and $W(.) = \frac{\varphi(.)}{\Phi(.)}$.

Solving simultaneously of the above equations is not mathematically tractable so one should apply some numerical optimization routine to get the solutions.

### 6.2. Real Life Applications

Here we have considered two datasets which are related to N latitude degrees in 69 samples from world lakes, which appear in Column 5 of the Diversity data set in website: http://users.stat.umn.edu/sandy/courses/8061/datasets/lakes.lsp and the white cells count (WCC) of 202 Australian athletes, given in Cook and Weisberg (1994) for the purpose of data fitting. Using GenSA package in R we have fitted our proposed distribution, i.e., $GABSN(\alpha,\beta,\lambda,\mu,\sigma)$ along with the normal $N(\mu,\sigma^2)$ distribution, the logistic $LG(\mu,\beta)$ distribution, the Laplace $La(\mu,\beta)$ distribution, the skew-normal $SN(\lambda,\mu,\sigma)$ distribution of Azzalini (1985), the skew-logistic $SLG(\lambda,\mu,\beta)$ distribution of Wahed and Ali (2001), the skew-Laplace $SLa(\lambda,\mu,\beta)$ distribution of Nekoukhou and Alamatsaz (2012), the alpha-skew-normal $ASN(\alpha,\mu,\sigma)$ distribution of Elal-Olivero (2010), the alpha-skew-Laplace $ASLa(\alpha,\mu,\beta)$ distribution of Harandi and Alamatsaz (2013), the alpha-skew-logistic $ASLG(\alpha,\mu,\beta)$ distribution of Hazarika and Chakraborty (2014), the alpha-beta-skew-normal $ABSN(\alpha,\beta,\mu,\sigma)$ distribution and beta-skew-normal $BSN(\beta,\mu,\sigma)$ distribution of Shafiei et al. (2016), and the generalized alpha-skew-normal $GASN(\alpha,\lambda,\mu,\sigma)$ distribution of Sharafi et al. (2017) for comparison purpose.

The values of the MLE's of the parameters for different distributions along with log-likelihood, AIC and BIC are given in Tables 1 and 2.

**Table 1:** MLE's, log-likelihood, AIC and BIC for N latitude degrees in 69 samples from world lakes.

| Parameters → Distribution ↓ | $\mu$ | $\sigma$ | $\lambda$ | $\alpha$ | $\beta$ | $\log L$ | AIC | BIC |
|---|---|---|---|---|---|---|---|---|
| $N(\mu,\sigma^2)$ | 45.165 | 9.549 | -- | -- | -- | -253.599 | 511.198 | 515.666 |
| $LG(\mu,\beta)$ | 43.639 | -- | -- | -- | 4.493 | -246.65 | 497.290 | 501.758 |
| $SN(\lambda,\mu,\sigma)$ | 35.344 | 13.7 | 3.687 | -- | -- | -243.04 | 492.072 | 498.774 |
| $BSN(\beta,\mu,\sigma)$ | 54.47 | 5.52 | -- | -- | 0.74 | -242.53 | 491.060 | 497.760 |
| $SLG(\lambda,\mu,\beta)$ | 36.787 | -- | 2.828 | -- | 6.417 | -239.05 | 490.808 | 490.808 |
| $La(\mu,\beta)$ | 43.0 | -- | -- | -- | 5.895 | -239.25 | 482.496 | 486.964 |
| $ASLG(\alpha,\mu,\beta)$ | 49.087 | -- | -- | 0.861 | 3.449 | -237.35 | 480.702 | 487.404 |
| $SLa(\lambda,\mu,\beta)$ | 42.3 | -- | 0.255 | -- | 5.943 | -236.90 | 479.799 | 486.501 |
| $ASLa(\alpha,\mu,\beta)$ | 42.3 | -- | -- | -0.22 | 5.44 | -236.08 | 478.159 | 484.861 |
| $ASN(\alpha,\mu,\sigma)$ | 52.147 | 7.714 | -- | 2.042 | -- | -235.37 | 476.739 | 483.441 |
| $ABSN(\alpha,\beta,\mu,\sigma)$ | 53.28 | 9.772 | -- | 2.943 | -0.292 | -234.36 | 476.719 | 485.655 |
| $GASN(\alpha,\lambda,\mu,\sigma)$ | 56.319 | 8.544 | -0.672 | 12.052 | -- | -230.531 | 469.062 | 477.998 |
| $GABSN(\alpha,\beta,\lambda,\mu,\sigma)$ | 58.333 | 6.489 | -0.616 | -6.144 | -5.870 | -225.369 | **460.738** | **471.909** |

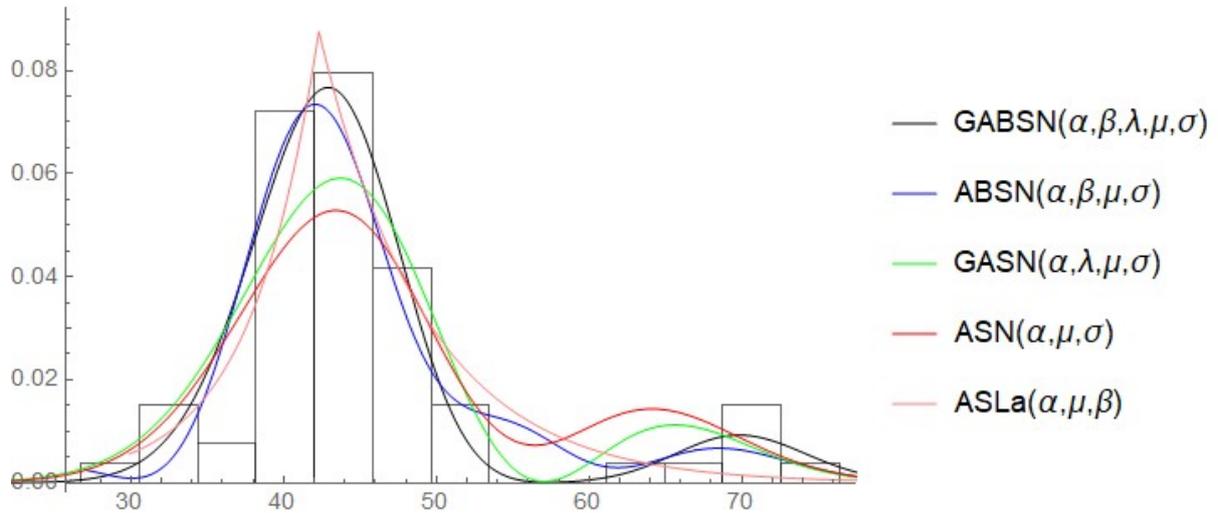

**Figure 3:** Plots of observed and expected densities of some distributions for N latitude degrees in 69 samples from world lakes.

**Table 2:** MLE's, log-likelihood, AIC and BIC for white cells count (WCC) of 202 Australian athletes.

| Parameters → <br> Distribution ↓ | $\mu$ | $\sigma$ | $\lambda$ | $\alpha$ | $\beta$ | $\log L$ | AIC | BIC |
|---|---|---|---|---|---|---|---|---|
| $N(\mu,\sigma^2)$ | 7.109 | 1.796 | -- | -- | -- | -404.919 | 813.838 | 820.455 |
| $La(\mu,\beta)$ | 6.844 | -- | -- | -- | 1.38 | -407.142 | 818.284 | 824.901 |
| $ASLa(\alpha,\mu,\beta)$ | 6.4 | -- | -- | -0.265 | 1.263 | -400.992 | 807.984 | 817.909 |
| $LG(\mu,\beta)$ | 6.996 | -- | -- | -- | 0.991 | -401.612 | 807.224 | 813.841 |
| $SLa(\lambda,\mu,\beta)$ | 6.4 | -- | 0.762 | -- | 1.287 | -400.366 | 806.732 | 816.657 |
| $BSN(\beta,\mu,\sigma)$ | 6.813 | 1.687 | -- | -- | -0.061 | -399.638 | 805.276 | 815.201 |
| $ASLG(\alpha,\mu,\beta)$ | 6.413 | -- | -- | -0.21 | 0.949 | -398.943 | 803.886 | 813.811 |
| $ASN(\alpha,\mu,\sigma)$ | 8.195 | 1.684 | -- | 0.874 | -- | -398.393 | 802.786 | 812.711 |
| $GASN(\alpha,\lambda,\mu,\sigma)$ | 5.569 | 2.689 | 2.242 | 0.5115 | -- | -395.545 | 799.090 | 812.323 |
| $SN(\lambda,\mu,\sigma)$ | 5.105 | 2.691 | 2.729 | -- | -- | -396.161 | 798.322 | 808.247 |
| $SLG(\lambda,\mu,\beta)$ | 5.512 | -- | 1.736 | -- | 1.331 | -395.897 | 797.794 | 807.719 |
| $ABSN(\alpha,\beta,\mu,\sigma)$ | 7.758 | 1.796 | -- | 0.813 | -0.128 | -394.833 | 797.666 | 810.899 |
| $GABSN(\alpha,\beta,\lambda,\mu,\sigma)$ | 9.729 | 1.508 | -0.669 | 0.411 | 0.428 | -392.788 | 795.576 | 812.117 |

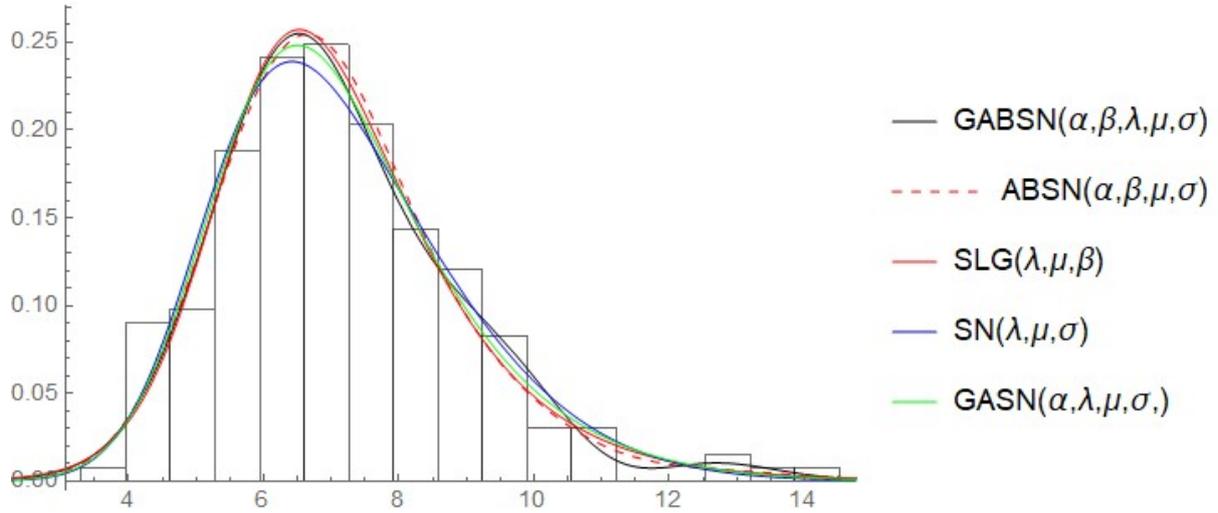

**Figure 4:** Plots of observed and expected densities of some distributions for white cells count (WCC) of 202 Australian athletes.

From Tables 1 and 2, it is seen that the proposed generalized alpha beta skew normal $GABSN(\alpha,\beta,\lambda,\mu,\sigma)$ distribution provides better fit to the data set under consideration in terms of all criteria, namely the log-likelihood, the AIC as well as the BIC. The plots of observed (in histogram) and expected densities (lines) presented in Figure 3 and 4, also confirms our finding.

### 6.2.1. Likelihood Ratio test:

Since $N(\mu,\sigma^2), SN(\lambda,\mu,\sigma)$, $ABSN(\alpha,\beta,\mu,\sigma)$, $GASN(\alpha,\lambda,\mu,\sigma)$, $GBSN(\beta,\lambda,\mu,\sigma)$ and $GABSN(\alpha,\beta,\lambda,\mu,\sigma)$ distributions are nested models, the likelihood ratio (LR) test is used to discriminate between them. The LR test is carried out to test the following hypothesis as shown in Table 3 below.

**Table 3:** The values of LR test statistic for different hypothesis.

| Hypothesis | LR values | |
|---|---|---|
| | Dataset 1 | Dataset 2 |
| $H_0 : \beta = 0$ vs $H_1 : \beta \neq 0$ | 10.324 | 5.514 |
| $H_0 : \lambda = 0$ vs $H_1 : \lambda \neq 0$ | 17.982 | 4.09 |
| $H_0 : \alpha = \beta = 0$ vs $H_1 : \alpha = \beta \neq 0$ | 35.342 | 6.746 |
| $H_0 : \alpha = \beta = \lambda = 0$ vs $H_1 : \alpha = \beta = \lambda \neq 0$ | 56.460 | 24.262 |

Since all the values of LR test statistic for different hypothesis exceed the 95% critical value that is, 3.841. Thus there is evidence in favors of the alternative hypothesis that the sampled data comes from $GABSN(\alpha, \beta, \lambda, \mu, \sigma)$ distribution.

## 7. Conclusions

In this article the generalized-alpha-beta-skew-normal distribution which includes at most four modes is introduced and some of its basic properties are investigated. Some extensions of the proposed distribution along with some of their properties are discussed.

The numerical results of the modelling of two real life data set considered here has shown that the proposed $GABSN(\alpha, \beta, \lambda, \mu, \sigma)$ distribution provides better fitting in comparison to the other known distributions.